\documentstyle[12pt]{article}
\title{Ordinal Computers}
\author{Ryan Bissell-Siders\\Math Department Princeton\\
Princeton NJ 08544\\
\\email:
rcsiders@math.princeton.edu}
\begin{document}
\maketitle


\begin{abstract} Can a computer which runs for time $\omega^2$ compute
more than one which runs for time $\omega$?  No.  Not, at least, for the
infinite computer we describe.  Our computer gets more powerful when the
{\it set} of its steps gets larger.  We prove that they theory of second
order arithmetic cannot be decided by computers running to countable time.
\end{abstract}

\vspace{0.3cm} \noindent Section 1. Introduction; Undecidability of Arithmetic.

Our motivation is to build a computer that will store and manipulate 
surreal numbers. 
Hackelroad \cite{Hackleroad} and Lurie 
\cite{Lurie} examined at least two ways to 
compute surreals in finite time, and shown the difficulty of building
a field 
of surreals in which 
$x>y$, $x=y+z$, and $x=y \times z$
are decidable.
Likewise for reals.
In the recursive reals algebra
but not order is decidable.
And so it seems that
the question of whether the theory of either field under $+, \times ,<$ is
decidable, ought to refer to decidability by some class of computers that
can compute more than finite-time Turing Machines.  If we're going to talk
about whether a computer can decide {\it facts} about {\it numbers}, then
let's have a computer that can {\it construct} all the numbers we want to
talk about and {\it decide} the algebra and order relations.  
Computers running to time $\aleph_1$ can compute all reals,
and to ordinal time can compute all surreals.  Now, what facts about such
numbers can ordinal computers decide? 


We will prove, in sections 2 and 5, that polynomials with variables in the
integers and the reals cannot be decided by a computer running for
countable time.  This is a curious result, since polynomials over the
reals can be decided by Elimination of Quantifiers, and polynomials over
integers can be decided in countable time by simply checking all the
possible inputs.  Our theorem suggests that both methods are tight, or,
more precisely, that there is no way to join them together into a single
countable-time algorithm for polynomials in integers and reals.  To be
specific, we will show that it cannot be decided in countable time whether
there exists a real number x so that for all possible choices of some 20
integers, a polynomial in the integers is zero, and $g < n \times x < g+1$ -- an
inequality in the integers and the real.  

In section 3 we define a general notion of an ordinal computer, and in
section 4 we prove something about them: that what can be affirmed in time
$< i$ is equivalent to what can be defined by the sentences of order $i$ in a
language of arbitrarily-high order.  This is simply a generalization of
the idea that sentences involving existential quantification over the
reals ought to be affirmable by a computer which ran to an arbitrary
countable time, and deniable by a computer running to time $\aleph_1$.  To
stress this, we will write $\aleph_1$ throughout to identify both the first
set nonisomorphic to $\aleph_0$, though this is often written as the ordinal 
$\omega_1$. 

A reader interested in ordinal computers may read only the brief section
3.  A reader interested in our strange theorem may skip sections 3 and 4. 

We assume CH -- the cardinality of the reals is the same as the first
uncountable ordinal.  

%
%

Our strange theorem concerns the language $(R,Z,<,+, \times )$.  That is,
statements with variables ranging over $R$, with symbols $<,+, \times $ and a
predicate $Z$ which is true exactly on the integers.  As Alex Wilkie
pointed out to the author, this is just
the theory of second order arithmetic.
Let A be an 
algorithm which halts in countable time on just the false 
statements.  We intend to show A doesn't exist. 

To each algorithm A we can build a finite-time machine B which accepts
real number inputs (an infinite, pre-written tape).  On real inputs
$x$ and $r$, our machine B will run for time $\omega$ just in case 
$x$ encodes a run of A,
starting with input $r$ and halting after countable time.  If $x$ does not
code a run of A on $r$, then B will halt in finite time.  If A halts in
countable time, there is some $x$ on which B runs for time $\omega$.  So A 
halts in countable time just in case B doesn't halt in finite time. 

We prove in the last section that there exists a statement in the language
$(R,Z,+, \times ,<)$ which is true of $x,r$, just in case B runs for time
$\omega$ on inputs $x,r$.  The variables $x,r$ are just the free variables
of this statement.  The statement also contains some quantified integer
variables, but no quantified real variables.  Now A halts on $r$
just in case for some real $x$, B doesn't halt on $x,r$.  That is, there 
is some $x$ so that our statement (call it $\phi$) is true.  So A halts on 
statement $r$ just in case there exists a real $x$ so that $\phi$ is true of 
$x,r$.  

We will prove that there is a halting problem for computers
running to countable times.  In section 3 we show that there is no
computer which halts in countable time 
just in case its input corresponds to a
computer that does not halt in countable time.  But if we can determine in
countable time whether any statement of our language is true or not, then by
the equivalence shown in the last paragraph, we could determine which
computers halt in countable time.


This demonstrates a class of simple formulas of second-order arithmetic not
decidable in countable time.  This has corrolaries that can be stated
without reference to ordinal computers.  For instance, we prove that
the theory of second-order arithmetic is not model complete (assuming CH).
If it is model complete, then any formula is
equivalent to an existential formula.  
Any existential formula of second-order arithmetic
can be put in the form ``for some integer values, $p$'' where $p$ is a 
formula of $(R,<,+, \times )$.  We can decide this statement in countable
time by checking whether $p$ is true at any particular integer values, using
elimination of quantifiers for $(R,<,+, \times )$.  
Unfortunately, we never expected $(R,Z,<,+, \times )$ to be model complete.
The formula
$\forall k \exists p,q |x-p/q| < 1/kq^2$ defines the reals whose
continued fraction terms are unbounded.  It seems unlikely that
the complement of this set is existentially definable.  

\vspace{0.3cm} \noindent Section 2.  Reducing a countable -- time machine to a 
finite -- time machine.

Let a program A have finitely many instructions and keep ordinal
variables.  Each instruction may increment a variable, switch control
as two variables are equal or not, or stop the program.  
That is: ``$x++$,'' ``if $x=y$ goto l,'' or ``stop.''
At a limit time-ordinal, control returns to the $0$th command.  
At a limit time-ordinal, the value of each variable becomes the 
limit of the values that it has achieved.  

We will construct a program B which accepts a real variable $x$ iff it 
codes the run of A.  As a bit string, $x$ is a sequences of $1$'s, the 
number of $1$'s indicating a number, with $0$'s separating numbers.
B separates $x$ into three or more sequences.  
The first, $z$, encodes a map from $\omega$ to the timesteps of A. 
The second, which we will call $c$ for control, 
is a sequence of numbers corresponding to lines of the program A:
$c_i$ is the command that was active at time $i$.  
For each variable $x$ that A uses, $x_i$
is the value of the variable $x$ at time $i$.  
How can $z$, a list of finite numbers, 
code a map from $\omega$ to an infinite countable ordinal?  It is
actually a list of a statements written in a language that B can interpret
so that B accepts only those $z$ which code a map from $\omega$ to a
countable ordinal.  The statements of $z$ are: ``$n<m$'', ``$m$ is a limit 
ordinal,'' ``$m=n+1$,'' which 
occurs for each $n$ unless $n$ is the final element, in which case 
$z$ contains the fact:
``$n$ is the final element''.  In all of these statements $n$ and $m$
are finite numbers; $z$ codes a re-ordering of the finite numbers so that
they have the same order structure as the timesteps of A. 
The statement ``$m$ is the final element'' must appear first.  In this
way, B can check whether or not there is a final element.  Because $z$
contains explicit successor and limit statements, B can affirm, in a 
finite amount of time, that $n$ is a limit or that $n$ succeeds $m$.
We require that all
statements involving numbers less than $k$ occur before time $2k^2$.  
There are at most 2 statements about any 
particular $m$ and $n$, so there is some $z$ listing 
all statements about numbers less than $k$ before time $2k^2$.

When B learns that $m$ is the final element, it checks 
that $c_m$ is the stop command.  When B learns $n<m$, 
it checks that $c_n$ is not the stop command.  When B learns $m<n$ 
it checks that it hasn't already
learned $n<m$.  This insures $z$ is a partial ordering.  B checks that 
$m<n$ or $n<m$ occurs before $2(n+m)^2$.  This assures that $z$ is a 
total order.  When B learns that $m=n+1$ it checks that there is no
$l$ between $n$ and $m$.  This implies that $z$ is discrete.
When B learns that $m=n+1$ it checks that $n$ is less than $m$.
That is, the indices of $m$ and $n$
are in the same order as the values they encode.
This all implies that $z$ represents a discrete, wellordered total
order.
When B learns that $m=n+1$, it checks that $c_m$ is 
the correct instruction to follow $c_n$ and
that $x_m$ is derived from $x_n$ by applying rule $c_n$.  When it learns 
that $m$ is a limit ordinal, it checks that $c_m$ is $0$,
and that $x_m$ is the limit of $x_n$ for $n<m$.  But how can B check 
that the variables limit properly?

In
order to check that all variables limit to their appropriate values, 
B accepts two reals, $x$ and $x'$, for each variable $x$ used in 
A.  $x$
is, like all of our variables, a sequence of numbers, represented by
a string of 1's, separated by zeroes.  The $i$th number of $x$, $x_i$,
represents the value of the variable $x$ at the countable-ordinal time $z_i$. 
But this value may be infinite!  So $x_i$ really is the $z$ encoding of 
the value of $x$ at the $z$ encoding of time $i$.  B wants to check that if
time $i$ limits to time $j$, then $x_i$ limits to $x_j$.  This seems very
difficult, because in finite time B has no way of knowing that any
particular sequence of ordinal numbers limits to another ordinal.  Indeed,
B cannot even determine what any of the infinite numbers encoded by z are,
in finite time.  So to check that $x$ is continuous, B checks that $x$ is
monotone, and that $x$ and $x'$ are inverses.  Monotonicity means that 
if $m<n$ then $x_m \leq x_n$.  That $x$ has $x'$ as inverse means: 
if $x_m = n$, then $x'_n \leq m$; if $x'_m = n$, then $x_n = x_m$.  
The sequence $x$ is not strictly increasing,
and it will happen that $x_m$ is the same value for many consecutive
timesteps; this introduces the asymmetry between $x$ and $x'$. 

B checks that $c_i$, the string encoding which command is active, timesteps
appropriately, by checking that if $c_i$ is active, then the next string to
be active is $c_{i+1}$, or, if $c_i$ is a switch on a variable value, B finds
this variable value and checks whether $c_{i+1}$ or the alternate command
was active next.  At any limit ordinal, be checks that the zero-th command
was active.  B checks that $x_i$ behave correctly, as well, by checking that
$x_{i+1}$ is $x_i$ unless the command active at time $i$ is the command
"increment $x$", in which case, $x_{i+1} = x_i + 1$.

\vspace{0.3cm} \noindent Section 3.  Ordinal computers defined.

An ordinal computer runs for ordinal time, accepts ordinal inputs, and keeps
ordinal variables.  It has finitely many instructions of the form
``increment $x$'' or ``if $x=y$ goto instruction l'' or ``stop''
Minsky \cite{Minsky} proves that these are sufficient to compute 
all Turing Machines 
running to finite time.
Actually, he proves that ``increment $x$'', ``if $x=0$ goto l,'' and 
``decrement $x$'' are
sufficient.  But we can model ``decrement $x$'' with our more general goto
switch in a subroutine that starts with variables $a$ and $b$ equal to $0$.  
Variable $a$ is incremented. 
Then $a$ and $b$ are incremented until $a=y$.  Variable $b$ is returned; 
$b$ is the decrement
of $a$.  Our decrement subroutine, on an input without a predecessor, is the
identity.  But I don't think we can model the generalized goto switch
using decrement, increment, and ``if $x=0$ goto l''.  We can model the command
$x:= y$ by incrementing $x$ until it equals $y$.]

At a limit ordinal, what happens to the internal state of the machine? 
Command returns to instruction $0$.  Variables are set equal to their limit,
if they have one; otherwise they are set to zero. 

The halting problem is as difficult for computers halting at infinite
ordinal as it is for computers halting at finite ordinals.  Consider the
set of computers which {\it halt} when given themselves as input.  Let A
be a computer halting on exactly those computers which don't halt on
themselves.  Then run A on input A.  It halts iff it doesn't.  This is
true if we take halting to mean halting in finite time, ordinal time,
halting before $10$ timesteps have gone by, or before an uncountable number
of timesteps have passed.  In the final 
section, we will show that a computer A
halts in countable time just in case some statement of $(R,Z,+, \times ,<)$ is
true.  This implies that that language is undecidable in countable time
because to decide it would be to solve the countable-time halting problem
in countable time. 

\vspace{0.3cm} \noindent 
Section 4: The computational power of a run depends on the set of its
timesteps. 

We will prove that all ordinals which are equivalent under re-ordering
have the same computatinal power.  Let $\aleph_n$ be the first ordinal larger
as a set than $\aleph_i$ for $i<n$.  We want to know if algorithm A halts
before time $\aleph_i$ on input $a<\aleph_i$.  
There is a computer B which halts
before time $\aleph_j$, for some $j<i$, on input $a$ and
all inputs $b < \aleph_i$, just
in case A runs to time $\aleph_i$.  This will all be simpler if we set i=1. 
Then: A stops in countable time just in case B, on all real inputs $b$, does
NOT halt in finite time.  

B checks that $b$ codes a complete run of A.  That is, $b$ is a bit string
which encodes: 1. a map from some $\aleph_j$, for $j<i$, to the steps of A.
Program B sees the steps of A streaming by, in an order rearranged to be
as short as possible.
2. Which instruction of A was in command at each time.
3. The value of all the variables of A at all times.
4. An inverse for each variable, which encodes when the 
variable was $<$, $=$, or $>$ than each possible value.

The hard part to check is that the value of the variable at a limit time
$\lim_i$ is the limit of the values at times $i$.  We have already described
how such a B can operate, in the last paragraph of the previous section:
by checking that the variable $b$ encodes strings $x$ and $x'$ for each 
variable, which represent
inverses, and so that $x$ is monotone and invertible, hence continuous.  
We described the computer in the
previous section in great detail, and here it is all the same, but with
``less than $\aleph_i$'' replacing ``countable'' and ``less than some 
$\aleph_j$ for $j<i$'' replacing ``finite.''  

Our description of B in terms of A is entirely primitive recursive and not
dependent on A, so A may be considered a variable.  Indeed, there is a
primitive recursive algorithm to produce $A_n$ from $A_{n+1}$ so that 
$A_{n+1}$ halts on
input $a$ before time $\aleph_{n+1}$ iff for some $b$, $A_n$ 
doesn't halt before
time $\aleph_n$.  This allows us to describe the set of inputs on which a
program halts in terms of a arbitrarily-high-order language.  We start at
the level of predicates on finite numbers. 

\begin{center}
	$P(x_0)$ is true.	\\
	$A_0$ stops before $\aleph_0$ iff $\exists x_0 < \aleph_0 P(x_0)$ \\
	$A_1$ stops before $\aleph_1$ iff $\exists x_1 < \aleph_1
			\forall  x_0 < \aleph_0 P(x_0)$ \\
	$A_2$ stops before $\aleph_2$ iff $\exists x_2 < \aleph_2 
	      \forall x_1 < \aleph_1 \exists x_0 < \aleph_0, P(x_0)$ \\

	$A_\omega$ halts before $\aleph_\omega$ iff some sentence of the form \\
	 $\forall x_7 \exists x_6 \forall x_5 \exists x_4 \dots A0$ is true \\
\end{center}

We have reduced the set of inputs on which some computer halts to the set
of $x_0$ for which some high-level statement is true, so that runs halting
before some cardinal time decide sets which are of the same level in the
hierarchy. 

We remark that in the hierarchy above, sentences need not be so long. We
can exchange the existential and universal quantifiers if we are willing
to quantify over longer ordinals.  For instance: 
$\forall x<X \exists y<Y A0(x,y)$ is equivalent to
$\exists f<Y^X \forall x<X A0(x,f(x))$.  
So, $A_\omega$ halts before $\aleph_\omega$ iff some sentence of the form
$\exists x < \aleph_\omega \forall y < \aleph_\omega, A0(x,y)$ is true.


\vspace{0.3cm}
\noindent 
Section 5. Turning a finite - time machine which accepts real inputs into
a polynomial.

Let us remember that in section 1 we wanted to build a computer B out of a
computer A so that: computer A stops at countable time
on countably-long bit-string $r$ just in case 
there exists a countably-long bit-string
$x$ s.t. computer B doesn't stop in finite time on input $x,r$.  
The string $x$ codes a map from
omega to the timesteps A took; for each timestep of A, the value of each
variable, and which instruction was operating.  B keeps finite variables,
and is allowed to increment and compare them.  B is also allowed to switch
on the $i$th bit of its real input.
B was built in section 2.  Now we want
to code B with a polynomial relation of the form $p(Z)=0 \ \&\ q(R,Z)>0$, 
following Jones and Matijasevich \cite{Matijasevich}, so
that B halts iff its statement in the language of inequalities of
polynomials is true for no integers m.  That is:
A halts in countable time on $r$ iff
$\exists x$ B doesn't stop in finite time on $r,x$ iff
$\exists x \forall m \phi$.

That is, a computer running for countable time and keeping countable
variables will halt on those reals $r$ so that $\exists x \forall m \phi$. 
That will be proven once we turn B into a polynomial.  We turn to this
now. 

B is allowed the commands ``n++'' and ``if n=m goto'' for its finite 
variables
$n,m$.  It is also allowed to switch on the $i$-th bit of the real 
variable $x$:
``if $x_i=0$ goto'', where $i$ is a finite variable stored by B.  We want,
however, to weaken our program so it may only switch on the $i$-th bit 
of $x$
at time $n i$ for some integer $n$.  
This can be done by encoding B in an interface program.  The interface
is only allowed to switch on the time-indexed bit, but it successfully
stores all the bits of $x$, and then B can switch on a stored bit.  In
more detai:
We store $x$, as a
binary integer, until the program to store a bit takes longer than $n$
steps.  Then we execute B on the resulting integer, replacing ``if $x_i=0$
goto'' with the command to compute $2^i$ and bitwise multiply this by $x$, and
put the result in variable $y$.  Then ``if $y=0$ goto''.  
We watch that B stops
normally.  If B halts on $x$, then for some large enough $n$ B will halt
normally on the truncation of $x$.  If B doesn't halt on $x$, then for
no $n$ will B halt normally.
So: There exists $x$ so that our more powerful computer
stops on $x,r$ just in case there exists $x,n$ so that our weaker computer
stops on $x,r,n$.  Now change B to read off the alternate bits of a single
variable as $x$ and $n$.  So we have simplified B as desired.  

To this more restricted program B we associate a polynomial 
relation; a statement in the language $(R,Z,<, \times ,+)$.  B halts
in finite time on real number inputs $x,r$ just in case there are some
finite numbers $a,b,c,d \dots$ so that 
$p(B,r,x,a,b,c,d...) = 0 \ \&\ q(x,g) > 0$.  This is the form whose
instances we will prove form
an undecidable class of statements.
The first 
thing to do is to multiply $x$ and $r$ by $2^a$ so as to get numbers with 
positive integer parts (let's assume $x$ and $r$ have no positive 
part).  We will find $g=[2^a r]$ and $h=[2^a x]$ and henceforth only deal 
with g and h:
$\exists g$ s.t. $g<2^a r<g+1 \ \&\
\exists h$ s.t. $h<2^a x<h+1$
where our integer exponentiation is, by Matijasevich's famous proof, 
expressible as a polynomial relation. 

The rest of the polynomial can be interpreted as checking that $a,b,c,d 
\dots$
record a run of B which stops on ``inputs'' $g$ and $h$.  
If B stopped on inputs $x$ and $r$, then B will also stop
on some truncation of $x$ and $r$.  The integers
$b,c,d \dots$ are bit strings.  Substrings of length 
$n$ represent the state of
B.  So the first thing to do is stretch $g$ and $h$ out so that their 
bits are separted by $n-1$ 0's.  
$\exists i$ s.t. $i$ is the stretching of $g$ by a factor of $n$. 
$\exists j$ s.t. $j$ is the stretching of $h$ by a factor of $n$. 

Let us immediately prove that these can be coded by polynomials: We only
need prove that stretching of finite bit strings can be computed by a
Turing Maching; then it can be defined by polynomials. We take from
Minsky's paper the result that there is a Turing Maching which turns a bit
string $bbbbb$ into $b0b0b0b0b0$.  How?  Erase the leftmost one, 
and write it at
the same location on a second tape.  In this case, that means to turn
$1bbbb$ into $bbbb$, and write $10000$.  
The number on the second tape is $2^n$. 
Minsky's maching $W(2,3)$ turns this into $3^n$, which $W(3,4)$ 
turns into $2^2n$. 
Then we write this on a third tape, and start over.  When we are done,
$bbbbb$ has become $b0b0b0b0b0$.  

From here on we will follow Matijasevich and Jones.  
We need only add a
statement to take care of the commands ``1: if $x_i=0$, goto l''.  That is,
command l is active only if command 1 was previously active and $x$ is zero. 
But now since $x$ is finite and properly spaced, our command takes the form
of Matijasevich and Jones:
``Command l is bitwise dominated by Command 1 minus $x$''

So B halts iff $\exists a \dots  \exists j$, $a$ through $j$ all integers,
s.t. $p(B,a,b,c, \dots j,n) = 0 \ \& \ g<2^a x<g+1 \ \& \ h<2^a r < h+1$.  
Let $\phi$ be ``$p(B,a,b,c, \dots j,n) = 0 \ \&\ g<2^n x<g+1 \ \& \ 
 h<2^n r < h+1$.''
So A halts in countable time iff $\exists$ real $x$ s.t. $\forall$ 
integers $a \dots j \neg \phi$.

\end{document}